\def \NN {\mathbb N}
\def \CC {\mathbb C}

\def \RR {\mathbb R}
\def \ZZ {\mathbb Z}
\def \HH {\mathbb H}

\def \epsilon{\varepsilon}

\def \S  {{\mathcal S}}

\def \d {\text{d}}

\def \ep {\epsilon}
\def \ga {\gamma}
\def \Ga {\Gamma}
\def \si {\sigma}
\def \la {\lambda}

\renewcommand{\S}{{\mathcal S}}

\documentclass[12pt,reqno]{amsart}

\usepackage{amsmath,amssymb}

\usepackage{url} 

\numberwithin{equation}{section}

\begin{document}

%\date{preliminary version \today}

\title[]{$L$-functions of degree $2$ and conductor $1$: \\ underlying ideas and a generalisation} 

\author[]{J.KACZOROWSKI \lowercase{and} A.PERELLI}
%\date{}
\maketitle

%\bigskip
{\bf Abstract.} We present a streamlined account of a recent theorem on the classification of the $L$-functions of degree 2 and conductor 1 from the extended Selberg class. We also present a more general new result dealing with functional equations involving two Dirichlet series. Further, we correct a slip in the original proof of the above theorem, which however does not affect the final result.

%\medskip

%\medskip

%\bigskip
%\centerline{Contents}

\vskip.5cm
%-1-%%%%%%%%%%%%%%%%%%%%%%%%%%%%%%%%%%%%%%%%%%%%%%%%%%%%%%%%%%%%%%%%%%%%%%%%%%%%%%%%%%%%%%%%%%%%
\section{Introduction}

\smallskip
%%%%%%
The Selberg classes $\S$ and $\S^\sharp$ provide an excellent general framework for the study of the structure of $L$-functions from the analytic point of view, and several results have been obtained in this direction. We refer to our survey papers \cite{Kac/2006},\cite{Ka-Pe/1999b},\cite{Ka-Pe/2022a},\cite{Per/2005},\cite{Per/2004},\cite{Per/2010},\cite{Per/2017} for a rather exhaustive account of such results, and to the next section for definitions and some prerequisites. Of particular interest is the theory of nonlinear twists of the functions $F\in\S^\sharp$, which offers a great help in many situations. Applications of such a theory lead, for example, to a fully explicit description of the members of $\S^\sharp$ of degree $d<2$. Here we present a streamlined account of a recent result \cite{Ka-Pe/2022b} on the classification of the $L$-functions of degree 2 and conductor 1, highlighting the role played by the nonlinear twists. 

\smallskip
%%%%%
It is widely expected that the Selberg class $\S$ coincides, loosely speaking, with the class of automorphic $L$-functions; the results so far obtained for small degrees confirm this expectation. Much more mysterious is the content of the extended Selberg class $\S^\sharp$, for which, as far as we know, there is no conjectural description even in the case of degree 2. However, in the special case of conductor $q=1$ a full description is possible, thanks to a peculiar property of certain linear twists which holds in such a case. Quite possibly, an analogous description can be obtained for other small integer conductors, like $q=2$, 3 and 4.

\smallskip
In order to state our results we need to introduce the following normalization of the functions in $\S^\sharp$ of degree 2 and conductor 1. We say that $F\in\S^\sharp$ is {\it normalized} if its internal shift $\theta_F$ vanishes and the first nonvanishing Dirichlet coefficient equals 1. Our aim is to show that the functions $F\in\S^\sharp$ under inspection coincide with the (suitably normalized) $L$-functions of Hecke or Maass forms of level 1, whose functional equation reflects $F(s)$ to $F(1-s)$, contrary to the Selberg class case where $F(s)$ reflects to $\overline{F}(1-s)$. As we shall see, the normalized $F\in\S^\sharp$ with $d=2$ and $q=1$ have real coefficients, thus $\overline{F}=F$ and hence the two functional equations agree in this case. Moreover, every such $F\in\S^\sharp$ can be normalized, meaning that the normalization hypothesis does not affect the generality of our result. Indeed, if $F(s)$ is entire we shift $s$ by $-i\theta_F$ and divide by the first non-zero coefficient, while if $F$ has a pole at $s=1$ then already $\theta_F=0$ so we simply divide as before; see Lemma 4.1 of \cite{Ka-Pe/2022b} for these results.

\smallskip
%%%%%
The classification of the normalized functions $F\in\S^\sharp$ of degree 2 and conductor 1 is based on the value of the numerical invariant $\chi_F$ defined as
\begin{equation}
\label{chieffe}
\chi_F = H_F(1)+H_F(2)+2/3.
\end{equation}
In view of the definition of the above $H$-invariants, see \eqref{2-4} and \eqref{2-5} below, the value of $\chi_F$ can easily be computed from the data of the functional equation of $F$. For example, if $f$ (resp. $u$) is a Hecke cusp form  (resp. a Maass form) of level 1 with first nonvanishing Fourier (resp. Fourier-Bessel) coefficient equal to 1, and $L(s,f)$ (resp. $L(s,u)$) denotes its $L$-function, then a simple computation shows that

$\bullet$ if $f$ has weight $k$ then $F(s)=L(s+\frac{k-1}{2},f)$ is a normalized function in $\S^\sharp$ with $d=2$ and $q=1$, and
\[
\chi_F= \frac{(k-1)^2}{2};
\]

$\bullet$ if $u$ has weight 0 and eigenvalue $1/4+\kappa^2$ then again $F(s)=L(s,u)$ is a normalized function in $\S^\sharp$ with $d=2$ and $q=1$, and
\[
\chi_F=-2\kappa^2.
\]
It is interesting that, viceversa, the value of $\chi_F$ detects the nature of any normalized $F\in\S^\sharp$ of degree 2 and conductor 1, as the following result shows. Sometimes, thanks to the normalization of $F$, $\chi_F$ detects even the function $F$ itself; for example
\[
\chi_F=0 \ \Longrightarrow \ F(s)= \zeta(s)^2 \quad \text{and} \quad  \chi_F = \frac{121}{2} \Longrightarrow F(s) = L(s +\frac{11}{2},\Delta),
\]
where $\Delta$ is the Ramanujan cusp form.

\medskip
%%%%%
{\bf Theorem 1.} {\sl Let $F\in\S^\sharp$ of degree $2$ and conductor $1$ be normalized. Then $\chi_F\in\RR$ and

\smallskip
(i) if $\chi_F>0$ then there exists a holomorphic cusp form $f$ of level $1$ and even integral weight $k= 1 + \sqrt{2\chi_F}$ such that $F(s) = L\big(s+\frac{k-1}{2},f\big)$;

\smallskip
(ii) if $\chi_F=0$ then $F(s) = \zeta(s)^2$;

\smallskip
(iii) if $\chi_F<0$ then there exists a Maass form $u$ of level $1$, weight $0$ and with eigenvalue $1/4+\kappa^2= (1-2\chi_F )/4$ such that $F(s) = L(s,u)$.}

\medskip
%%%%%
Theorem 1 can be regarded as a sharp and general form of the classical converse theorems of Hecke \cite{Hec/1983} and Maass \cite{Maa/1949}. Actually, these theorems are used at the very end of our proof, since our strategy is, in broad terms, to show that eventually the functions $F$ in Theorem 1 satisfy a functional equation of Hecke or Maass type. The proof of Theorem 1 is achieved in four steps, and in this paper we analyse the ideas involved in each step. Then, in Section 7, we correct a slip in the original proof given in \cite{Ka-Pe/2022b}.

\smallskip
%%%%%
From Theorem 1 we see that every function in $\S^\sharp$ with $d=2$ and $q=1$ coincides, after suitable normalizations, with a level 1 modular $L$-function. For a non-normalized version of Theorem 1 we refer to Theorem 2 below. Moreover, since the Euler products in $\S$ are linearly independent (see \cite{K-M-P/1999} and \cite{K-M-P/2006}), the functions in the Selberg class $\S$ with degree $2$ and conductor $1$ are shifts of the (normalized) $L$-functions of Hecke or Maass eigenforms of level 1. In other words, they are automorphic $L$-functions, as expected.

\smallskip
%%%%%
We finally recall the following old conjecture in \cite[p.103]{Ka-Pe/2002}: {\it conductor, root number and the $H$-invariants with $n\leq d$ characterize the functional equation of any $F\in\S^\sharp$ of degree $d$}. Since in case (iii) of Theorem 1 the parity $\epsilon$ of $u$ can be expressed in terms of the root number as
\begin{equation}
\label{intro-0}
\ep= \frac{1-\omega_F}{2} \in\{0,1\},
\end{equation}
in view of \eqref{intro-0} our conjecture follows from Theorem 1 for all normalized $F\in\S^\sharp$ with $d=2$ and $q=1$.

\medskip
%%%%%
Hecke \cite{Hec/1983} studied also the case of two Dirichlet series $F,G$, with suitable analytic properties, satisfying a functional equation of type
\begin{equation}
\label{intro-1}
\Gamma(s) F(s) = \omega \Gamma(1-s) G(1-s),
\end{equation}
and showed that \eqref{intro-1} holds if and only if the Fourier series $f,g$ associated with $F,G$ satisfy a certain modularity relation. We conclude the survey adding the following new result, dealing with functional equations involving two Dirichlet series $F$ and $G$. Let $F,G$ be non identically vanishing and absolutely convergent for $\si>1$, and $(s-1)^mF(s),(s-1)^mG(s)$ be entire functions of finite order for some integer $m\geq0$. Moreover, let $F$ and $G$ satisfy a functional equation of type
\begin{equation}
\label{intro-2}
\ga(s)F(s) = \omega\overline{\ga}(1-s) G(1-s),
\end{equation}
where $|\omega|=1$ and the $\ga$-factor $\ga(s)$ is as in \eqref{2-2} below. Finally, let $d_\ga,q_\ga,\theta_\ga$ be defined as in \eqref{2-3} below, and let $\chi_\ga$ be defined as in \eqref{chieffe} by means of the data of the shifted $\ga$-factor $\ga(s-i\theta_\ga)$.

\medskip
%%%%%
{\bf Theorem 2.} {\sl Let $F,G$ be as above with $d_\ga=2$ and $q_\ga=1$. Then $\chi_\ga\in\RR$ and

(i) if $\chi_F>0$ then there exists a holomorphic cusp form $f$ of level $1$ and even integral weight $k= 1 + \sqrt{2\chi_F}$ such that $F(s) = L\big(s+\frac{k-1}{2}+i\theta_\ga,f\big)$ and $G(s) = L\big(s+\frac{k-1}{2}+i\theta_\ga,i^{k}\overline{\omega} f\big)$;

(ii) if $\chi_\ga=0$ then there exists $\alpha\in\CC\setminus\{0\}$ such that $F(s) = \alpha\zeta(s)^2$ and $G(s) = \alpha\overline{\omega} \zeta(s)^2$;

(iii) if $\chi_F<0$ then there exists a Maass form $u$ of level $1$, weight $0$, with eigenvalue $1/4+\kappa^2= (1-2\chi_F )/4$ and parity $\ep$ given by \eqref{intro-0} such that $F(s) = L(s+i\theta_\ga,u)$ and $G(s) = L(s+i\theta_\ga,(-1)^\ep\overline{\omega}u)$.}

\medskip
%%%%%
Clearly, Theorem 2 is more general than Theorem 1, but actually the two theorems are essentially equivalent. Indeed, in Section 8 we show that Theorem 2 is a consequence of Theorem 1 together with some features of the normalization process described in Section 3 and the fact that eventually the $\ga$-factors in \eqref{intro-2} turn out to be virtual $\ga$-factors as described in Section 4.

\medskip
%%%%%%
{\bf Acknowledgements.} This research was partially supported by the Istituto Nazionale di Alta Matematica, by the MIUR grant PRIN-2017 {\sl ``Geometric, algebraic and analytic methods in arithmetic''} and by grant 2021/41/BST1/00241 {\sl ``Analytic methods in number theory''}  from the National Science Centre, Poland.

\medskip
%-2-%%%%%%%%%%%%%%%%%%%%%%%%%%%%%%%%%%%%%%%%%%%%%%%%%%%%%%%%%%%%%%%%%%%%%%%%%%%%%%%%%%%%%%%%%%%%
\section{Definitions and prerequisites}

\smallskip
%%%%%
Throughout the paper we write $s=\si+it$, $\overline{f}(s)$ for $\overline{f(\overline{s})}$ and $e(x)$ for $e^{2\pi ix}$. The extended Selberg class $\S^\sharp$ consists of non identically vanishing Dirichlet series 
\[
F(s)= \sum_{n=1}^\infty \frac{a(n)}{n^s},
\]
absolutely convergent for $\si>1$, such that $(s-1)^mF(s)$ is entire of finite order for some integer $m\geq0$, and satisfying a functional equation of type
\begin{equation}
\label{2-1}
F(s) \gamma(s) = \omega \overline{\gamma}(1-s) \overline{F}(1-s),
\end{equation}
where $|\omega|=1$ and the $\gamma$-factor
\begin{equation}
\label{2-2}
\gamma(s) = Q^s\prod_{j=1}^r\Gamma(\lambda_js+\mu_j) 
\end{equation}
has $Q>0$, $r\geq0$, $\lambda_j>0$ and $\Re(\mu_j)\geq0$. The Selberg class $\S$ is, roughly speaking, the subclass of $\S^\sharp$ of the functions with Euler product and satisfying the Ramanujan conjecture $a(n)\ll n^\ep$. We refer to our survey papers \cite{Kac/2006},\cite{Ka-Pe/1999b},\cite{Ka-Pe/2022a},\cite{Per/2005},\cite{Per/2004},\cite{Per/2010},\cite{Per/2017} for further definitions, examples and the basic theory of the Selberg class. 

\smallskip
%%%%%%
Degree $d$, conductor $q$, root number $\omega_F$ and $\xi$-invariant $\xi_F$ of $F\in\S^\sharp$ are defined by
\begin{equation}
\label{2-3}
\begin{split}
d=d_F :=2\sum_{j=1}^r\lambda_j, \qquad q=q_F:= (2\pi)^dQ^2\prod_{j=1}^r\lambda_j^{2\lambda_j}, \\
\omega_F=\omega \prod_{j=1}^r \lambda_j^{-2i\Im(\mu_j)},  \qquad \xi_F = 2\sum_{j=1}^r(\mu_j-1/2):= \eta_F+ id\theta_F
\end{split}
\end{equation}
with $\eta_F,\theta_F\in\RR$; $\theta_F$ is the internal shift of $F$. The $H$-invariants are defined for every $n\geq0$ as
\begin{equation}
\label{2-4}
H_F(n) = 2 \sum_{j=1}^r \frac{B_n(\mu_j)}{\lambda_j^{n-1}},
\end{equation}
where $B_n(z)$ is the $n$-th Bernoulli polynomial. Clearly
\begin{equation}
\label{2-5}
H_F(0)=d, \qquad H_F(1)=\xi_F,  \qquad H_F(2) = 2 \sum_{j=1}^r \frac{\mu_j^2 -\mu_j + 1/6}{\lambda_j}.
\end{equation}

\smallskip
%%%%%
In the rest of this section we restrict ourselves to normalized $F\in\S^\sharp$ with $d=2$ and $q=1$. Indeed, although several statements below hold in greater generality, in the case at hand they are simpler and more explicit; proofs can be found in \cite{Ka-Pe/2005},\cite{Ka-Pe/2017},\cite{Ka-Pe/2021}.

\smallskip
%%%%%
We already remarked in the Introduction that any normalized $F$ has real coefficients, so $\overline{F}=F$ and its $H$-invariants are real. Functional equation \eqref{2-1} can be written in the asymmetric form
\begin{equation}
\label{invform}
F(s) = \omega S_F(s)h_F(s) F(1-s),
\end{equation}
where
\begin{equation}
\label{5-14}
S_F(s) := 2^r \prod_{j=1}^r \sin(\pi(\lambda_js + \mu_j)) = \sum_{j=0}^N a_j e^{i\pi\omega_js}
\end{equation}
with
\begin{equation}
\label{5-16}
N\geq 1, \quad -1=\omega_0<\omega_1<\cdots<\omega_N=1, \quad \omega_j=-\omega_{N-j}, \quad a_j\neq0, \quad a_0=a_N=-\omega_F,
\end{equation}
and
\begin{equation}
\label{h-def}
h_F(s) = \frac{1}{(2\pi)^r} Q^{1-2s} \prod_{j=1}^r \big(\Gamma(\lambda_j(1-s)+\overline{\mu}_j) \Gamma(1-\lambda_js-\mu_j)\big).
\end{equation}
Note the slight difference with the definition of the $h$-function $h_F(s)$ given in our previous papers, which corresponds to $\omega h_F(s)$. The reason for such a minor change is that \eqref{h-def} depends only on the $\gamma$-factor, without involving the $\omega$-datum, and this will turn out to be convenient in Section 4 when defining the $h$-function associated with virtual $\gamma$-factors. The functions $S_F(s)$ and $\omega h_F(s)$ are invariants, and $h_F(s)$ has the expansion
\begin{equation}
\label{h-asy}
h_F(s) \approx \frac{\omega_F}{\omega} \frac{(4\pi)^{2s-1}}{\sqrt{2\pi}} \sum_{\ell=0}^\infty d_F(\ell) \Gamma(3/2-2s-\ell)
\end{equation}
where $ \approx$ means that cutting the sum at $\ell=M$ one gets a meromorphic remainder which is $\ll$ than the modulus of the $M$-th term times $1/|s|$. The coefficients $d_F(\ell)$ are the {\it structural invariants} of $F$ and play a prominent role in the proof of Theorem 1; we have that $d_F(0)=1$.

\smallskip
%%%%%
The standard twist of $F$ is defined for $\si>1$ and $\alpha>0$ by
\begin{equation}
\label{standtwist}
F(s,\alpha) = \sum_{n=1}^\infty \frac{a(n)}{n^s} e(-\alpha \sqrt{n}),
\end{equation}
and the spectrum of $F$ is
\begin{equation}
\label{spec}
\rm{Spec}(F) :=\{\alpha>0: a(\alpha^2/4)\neq0\} = \{2\sqrt{m}: m\in\NN \ \text{with} \ a(m)\neq 0\}.
\end{equation}
Moreover, we write
\begin{equation}
\label{poles}
s_\ell = \frac{3}{4} - \frac{\ell}{2}, \quad \ell=0,1,2, \dots
\end{equation}
We recall that $F$ has polynomial growth on vertical strips. Moreover, the standard twist $F(s,\alpha)$ is entire if $\alpha\not\in$ Spec$(F)$, while for $\alpha\in$ Spec$(F)$ it is meromorphic on $\CC$ with at most simple poles at the points $s=s_\ell$, with residue
\begin{equation}
\label{resi}
 \rho_\ell(\alpha) =  \frac{e^{i\pi/4}\overline{a(\alpha^2/4)}}{\sqrt{\alpha}} d_F(\ell) \, (-2\pi i)^{-\ell} \alpha^{-\ell}.
\end{equation}
Note that $\rho_0(\alpha)\neq0$ for every $\alpha\in$ Spec$(F)$. Finally, $F(s,\alpha)$ has also polynomial growth on every vertical strip. We refer to \cite{Ka-Pe/2005},\cite{Ka-Pe/2016a},\cite{Ka-Pe/2021} for these and other results on $F(s,\alpha)$.

\medskip
%-3-%%%%%%%%%%%%%%%%%%%%%%%%%%%%%%%%%%%%%%%%%%%%%%%%%%%%%%%%%%%%%%%%%%%%%%%%%%%%%%%%%%%%%%%%%%%%
\section{Step 1:  nonlinear twists and invariants}

\smallskip
%%%%%
We start with the observation that, in view of \eqref{invform}, the function $h_F(s)$ contains a good amount of information on the functional equation of $F$. In turn, thanks to \eqref{h-asy}, $h_F(s)$ is determined by the structural invariants $d_F(\ell)$ up to a multiplicative constant; see the proof of Proposition 5.1 in \cite{Ka-Pe/2022b}. Moreover, the structural invariants appear explicitly in the residues \eqref{resi} of the standard twist $F(s,\alpha)$, see \eqref{standtwist}, an important member of the family of nonlinear twists. The nonlinear twists are a major source of information on the invariants of $\S^\sharp$, thanks to a transformation formula relating, loosely speaking, a nonlinear twist to its dual nonlinear twist; see \cite{Ka-Pe/2011},\cite{Ka-Pe/2016a},\cite{Ka-Pe/2016b}. Although in general such a transformation formula is quite complicated, in suitably chosen cases it is possible to derive from it some explicit information on the invariants. The reason lies in the fact that the transformation formula is highly sensitive to formal changes in the involved twists, thus producing somewhat unexpected and nontrivial identities. We illustrate this phenomenon in the special situation required by the proof of Theorem 1. 

\smallskip
%%%%%
For $\si>1$ and $\alpha,\beta\in\RR$ we consider the nonlinear twist
\[
F(s,\beta,\alpha) = \sum_{n=1}^\infty \frac{a(n)}{n^s} e(-\beta n-\alpha\sqrt{n}),
\]
where $a(n)$ are the Dirichlet coefficients of $F$. Observe that if $\beta\in\ZZ$ then $F(s,\beta,\alpha)$ reduces to the standard twist $F(s,\alpha)$ of $F$; in particular, $F(s,1,\alpha) = F(s,\alpha)$. Then, applying the transformation formula to $F(s,1,\alpha)$ with $\alpha\in$ Spec$(F)$, after long and painful computations to make explicit all the involved quantities, we get the following result. 

\smallskip
%%%%%
For any integer $m\geq0$ we define the polynomials in the $(s,\alpha)$-variables
\begin{equation}
\label{4-1}
\begin{split}
W_m(s,\alpha) &=\mathop{\sum_{\nu=0}^\infty \sum_{\mu=0}^\infty \sum_{k=3\nu}^\infty  \sum_{\ell=0}^{\infty} \sum_{h=0}^\infty}_{\substack{2|(\mu+k) \\ -2\nu+\mu+k+2\ell+h=m}}  A(\nu,\mu,k,\ell,h) \\
&\times {-\frac14-s-\frac{\ell}{2} \choose \mu}  {\frac12-2s +2\nu-\mu-k-\ell \choose h} d_F(\ell) \alpha^h,
 \end{split}
\end{equation}
where $d_F(\ell)$ are the structural invariants of $F$,
\[
A(\nu,\mu,k,\ell,h) = \frac{1}{\sqrt{\pi} \nu!} \big(-\frac{2i}{\pi}\big)^{\frac{\mu+k}{2}} \big(-\frac12\big)^h (4\pi)^{\nu-\ell} a_{k,\nu} \Gamma\big(\frac{\mu+k+1}{2}\big) i^{\nu+\ell}
\]
and the coefficients $a_{k,\nu}$ are defined by the expansion
\[
\Big(\sum_{k=3}^\infty {1/2 \choose k} \xi^k\Big)^\nu = \sum_{k=3\nu}^\infty a_{k,\nu} \xi^k.
\]
Note that the functions $W_m(s,\alpha)$ are polynomials since $k\geq 3\nu$, and hence equation $-2\nu+\mu+k+2\ell+h=m$ in \eqref{4-1} has only finitely many solutions for every $m$.

\medskip
{\bf Lemma 3.1.}  {\sl Let $F$ be as in Theorem $1$ and $\alpha\in$\hskip.1cm{\rm Spec}$(F)$. Then for every integer $M\geq0$ we have
\[
F(s,\alpha) =  \sum_{m=0}^{M} W_m(s,\alpha) F\big(s+\frac{m}{2},\alpha\big) + H_{M}(s,\alpha),
\]
where $W_0(s,\alpha)\equiv 1$ and the function $ H_{M}(s,\alpha)$ is holomorphic for $\si>-(M-1)/2$.}

\medskip
%%%%%
We do not discuss the proof of this result here, as a sketch of it is already presented in Lemma 4.2 of \cite{Ka-Pe/2022b}. 

\smallskip
%%%%%
The interest of Lemma 3.1 lies in the fact that the first term of the sum on the right hand side cancels with the left hand side, hence the function
\begin{equation}
\label{sumtwists}
 \sum_{m=1}^{M} W_m(s,\alpha) F\big(s+\frac{m}{2},\alpha\big) 
\end{equation}
is holomorphic for  $\si>-(M-1)/2$. But for $\alpha\in$ Spec$(F)$, see \eqref{spec}, $F(s,\alpha)$ has (at most) simple poles at the points $s=s_\ell$, see \eqref{poles}, with residues $\rho_\ell(\alpha)$ in \eqref{resi}. Since the residue of \eqref{sumtwists} at $s=s_M$ must be zero, recalling that $a(\alpha^2/4)\neq0$ we derive the identities
\begin{equation}
\label{4-27}
\sum_{m=1}^M W_m(s_M,\alpha) (-2\pi i)^m d_F(M-m) \alpha^m=0
\end{equation}
for every $\alpha\in$ Spec$(F)$ and for every integer $M\geq1$. Now we observe that the left hand side of \eqref{4-27} is a polynomial in $\alpha$, say $W_M(\alpha)$. Moreover, Spec$(F)$ is an infine set since $F$, having positive degree, cannot be a Dirichlet polynomial. Thus all coefficients of $W_M(\alpha)$ must vanish.

\smallskip
Next we note that the coefficients of $W_M(\alpha)$ are quadratic in the $d_F(\ell)$, since the polynomials $W_m(s,\alpha)$ in \eqref{4-1} are linear in $d_F(\ell)$. Actually, a computation shows that the $M$-th coefficient of $W_M(\alpha)$ can be arranged into the following quadratic form in the variables $X_\ell=d_F(\ell)$
\[
\widetilde{Q}_N(X_0,\dots,X_N) = \sum_{0 \leq \ell\leq N} \ \sum_{0\leq h \leq N-\ell}  (-2\pi i)^{-h} B_{2N-h}(s_{2N},\ell,h) X_\ell X_h,
\]
where $M=2N$ and $B_{2N-h}(s_{2N},\ell,h)$ are certain complicated but explicit coefficients. Moreover, the coefficients of $\widetilde{Q}_N(X_0,\dots,X_N)$ are independent of $F$ since so are the polynomials $W_m(s,\alpha)$, and it turns out that for $N\geq2$
\[
(-2\pi i)^{-N} B_N(s_{2N},0,N) + B_{2N}(s_{2N},N,0) \neq 0.
\]
Hence, after an obvious normalization of $\widetilde{Q}_N$ leading to a new quadratic form $Q_N$, the structural invariants of any $F$ as in Theorem 1 satisfy
\begin{equation}
\label{quadratic}
Q_N(d_F(0),\dots,d_F(N)) = \sum_{\substack{\ell,h\geq 0 \\ \ell+h\leq N}} \alpha_{\ell,h} d_F(\ell) d_F(h)=0,
\end{equation}
where the coefficients $\alpha_{\ell,h}$ are independent of $F$, can be shown to be real and satisfy $\alpha_{0,N}+\alpha_{N,0}=1$; see the proof of Proposition 4.1 in \cite{Ka-Pe/2022b} for all the above statements.

\smallskip
%%%%%
Equation \eqref{quadratic} is very interesting. Indeed, since $d_F(0)=1$, by induction we have that 
\begin{equation}
\label{induction}
\text{{\it all $d_F(\ell)$ with $\ell\geq2$ can be expressed in terms of $d_F(1)$, namely $d_F(\ell)= E_\ell(d_F(1))$,}}
\end{equation}
say, and such an expression $E_\ell(X)$ is independent of $F$. Moreover, \eqref{quadratic} shows that under the hypotheses of Theorem 1, the structural invariants $d_F(\ell)$ lie on a universal family of quadratic varieties. Quite possibly, the general case behaves similarly, namely the structural invariants of any $F\in\S^\sharp$ lie on a family of algebraic varieties depending on $F$ in a very mild way, perhaps only on few invariants of $F$. Since the structural invariants are closely connected with the $\gamma$-factor of $F$, this phenomenon, if true, could explain why the functional equations of $L$-functions are expected to involve only $\ga$-factors of very special shape.

\smallskip
%%%%%
Now we proceed with the computation of $d_F(1)$ in terms of the $H$-invariants, based on the Stirling expansion of the $\Ga$-functions in the definition \eqref{h-def} of $h_F(s)$. This is possible since both the Stirling expansion of $h_F(s)$ and the $H$-invariants involve the Bernoulli polynomials. Actually, the invariants $d_F(\ell)$ and $H_F(n)$ play a similar role in the Selberg class theory, in the sense that they are the coefficients of two asymptotic expansions related to the functional equation of $F$. The $H$-invariants carry more information and essentially determine the $\ga$-factor of $F$, while the structural invariants essentially determine only the $h$-function in the asymmetric form of the functional equation; see Theorems 1 and 2 in \cite{Ka-Pe/2022a}. After a careful computation and then comparing the output with formula \eqref{h-asy} for $h_F(s)$, we obtain that
\begin{equation}
\label{d_1}
d_F(1) = H_F(1)+H_F(2) + \frac{13}{24}.
\end{equation}
We remark that, in principle, this computation allows to determine explicitly all the $d_F(\ell)$.
Thus, in view of \eqref{chieffe} we get
\begin{equation}
\label{d-chi}
d_F(1)=\chi_F -\frac{1}{8}.
\end{equation}
Note that \eqref{d_1} explains the presence of the term $H_F(1)+H_F(2)$ in the definition \eqref{chieffe} of $\chi_F$, while the value $2/3$ in \eqref{chieffe} is a normalization suggested by \eqref{5-2} below, aiming at a nicely symmetric formulation of Theorem 1. In conclusion, from \eqref{induction} and \eqref{d-chi} the value of $\chi_F$ determines all the structural invariants and hence $\frac{\omega}{\omega_F}h_F(s)$. Moreover, since after normalization $F$ has real coefficients, the $H$-invariants are real, and so is $\chi_F$.

\smallskip
%%%%%
We summarise the output of the first step of the proof by the following proposition.

\medskip
%%%%%
{\bf Proposition 3.1.} {\sl Let $F$ be as in Theorem $1$. Then $\chi_F$ is real and its value determines the function $\frac{\omega}{\omega_F}h_F(s)$.}

\medskip
%-4-%%%%%%%%%%%%%%%%%%%%%%%%%%%%%%%%%%%%%%%%%%%%%%%%%%%%%%%%%%%%%%%%%%%%%%%%%%%%%%%%%%%%%%%%%%%%
\section{Step 2: virtual $\gamma$-factors}

\smallskip
%%%%%
As already mentioned in the Introduction, our ultimate goal is to show that every $F$ as in Theorem 1 satisfies a functional equation of Hecke or Maass type. Therefore, we introduce the {\it virtual $\ga$-factors}
\begin{equation}
\label{5-1}
\gamma(s) = 
\begin{cases}
(2\pi)^{-s} \Gamma(s+\mu) &\text{with $\mu>0$} \\
\pi^{-s} \Gamma\big(\frac{s+\epsilon+i\kappa}{2}\big) \Gamma\big(\frac{s+\epsilon-i\kappa}{2}\big) &\text{with $\epsilon\in\{0,1\}$ and $\kappa\geq0$,}
\end{cases}
\end{equation}
respectively of Hecke and Maass type. We denote the analog for such $\ga$-factors of the invariants of $F$ by adding the suffix $\ga$ instead of $F$. For example, we see that the values of degree and conductor of the $\ga(s)$ in \eqref{5-1} are 2 and 1, respectively, and a computation shows that
\begin{equation}
\label{5-2}
\quad\chi_\gamma:=  H_\gamma(1) + H_\gamma(2) +\frac23= 
\begin{cases}
2\mu^2 \\
- 2\kappa^2.
\end{cases}
\end{equation}
Moreover, the analog $h_\ga(s)$ of the function $h_F(s)$ in \eqref{h-def} has the expansion
\begin{equation}
\label{5-8}
h_\gamma(s) \approx \frac{(4\pi)^{2s-1}}{\sqrt{2\pi}} \sum_{\ell=0}^\infty d_\gamma(\ell) \Gamma(3/2 -2s-\ell), \quad d_\gamma(0)=1.
\end{equation}

\smallskip
We call virtual the $\ga$-factors in \eqref{5-1} since it is well known from the converse theorems of Hecke and Maass that not every such $\ga$-factor corresponds to an existing $L$-function. Nonetheless, their structural invariants $d_\ga(\ell)$ satisfy the same formal properties of the $d_F(\ell)$ in equations \eqref{induction}, \eqref{d-chi} and in Proposition 3.1, namely

\medskip
{\bf Lemma 4.1.} {\sl For every $\ell\geq2$ we have that $d_\ga(\ell) = E_\ell(d_\ga(1))$, $d_\ga(1)=\chi_\ga -1/8$ and the value of $\chi_\ga$ determines $h_\ga(s)$.}

\medskip
%%%%
The key point in the proof of Lemma 4.1 is the fact that the $d_\ga(\ell)$ are polynomials in the parameters $\mu$ or $\kappa$ in \eqref{5-1}, i.e there exist $P_\ell,Q_\ell\in\RR[x]$ such that
\begin{equation}
\label{5-9}
d_\ga(\ell) = 
\begin{cases}
P_\ell(\mu) \\
Q_\ell(\kappa),
\end{cases}
\end{equation}
the latter polynomials being independent of the value of $\ep$ in \eqref{5-1}. This is obtained by the same type of computation leading to \eqref{d_1}, which is based on the Stirling expansion and hence does not depend on the function $F$ and involves the Bernoulli polynomials. Next we note that $\mu$ and $\kappa$ in \eqref{5-1} are closely related to the weight of the Hecke forms and to the eigenvalue of the Maass forms, respectively, and that the sets of such weights and eigenvalues are infinite. Moreover, for these specific values of $\mu$ and $\kappa$, the associated normalized $L$-functions satisfy the hypotheses of Theorem 1, hence the corresponding $d_\ga(\ell)$ satisfy \eqref{quadratic}. But, in view of \eqref{5-9}, \eqref{quadratic} becomes a polynomial equation in the variable $\mu$ or $\kappa$, and has infinitely many solutions. Thus \eqref{quadratic} is satisfied for every $\mu$ and $\kappa$, and Lemma 4.1 follows by the same argument at the end of the previous section.

\smallskip
%%%%%
In view of \eqref{5-1} and \eqref{5-2}, the values of $\chi_\ga$ range over $\RR$, hence to every $F$ as in Theorem 1 we can associate a unique virtual $\ga$-factor $\ga(s)$ such that $\chi_\ga=\chi_F$. The uniqueness is obvious if $\chi_F>0$, while if $\chi_F\leq0$ we choose the virtual $\ga$-factor with the $\ep$ such that $\omega_F=(-1)^\ep$. This is possible since in \cite{Ka-Pe/2017} we proved, under the hypotheses of Theorem 1, that $\omega_F=\pm1$. Therefore, thanks to Proposition 3.1, Lemma 4.1, \eqref{h-asy} and \eqref{5-8} we have that
\begin{equation}
\label{equal-h}
h_F(s) = \frac{\omega_F}{\omega} h_\ga(s),
\end{equation}
$\ga(s)$ being the virtual $\ga$-factor associated with $F$. As a consequence, by \eqref{invform},\eqref{equal-h} and the definitions of $S_\ga(s)$ and $h_\ga(s)$ we can write
\[
\begin{split}
F(s) = \omega_F S_F(s)h_\ga(s) F(1-s) &= \omega_F \frac{S_F(s)}{S_\ga(s)} S_\ga(s)h_\ga(s) F(1-s) \\
&= \omega_F R(s) \frac{\ga(1-s)}{\ga(s)} F(1-s),
\end{split}
\]
where
\begin{equation}
\label{r-funct}
R(s) = \frac{S_F(s)}{S_\ga(s)}.
\end{equation}
Hence, if $F$ is as in Theorem $1$ and $\ga(s)$ is its associated virtual $\ga$-factor, then $F$ satisfies the functional equation
\begin{equation}
\label{virtual-fe}
\ga(s)F(s) = \omega_F R(s) \ga(1-s)F(1-s).
\end{equation}
This shows that the functional equation of $F$ is quite close to a functional equation of Hecke or Maass type, so we are approaching our goal. Actually, if the function $R(s)$ in \eqref{r-funct} is constant, then Theorem 1 follows immediately from the converse theorems of Hecke and Maass type already known in the literature; see Hecke \cite{Hec/1983} and Raghunathan \cite{Rag/2010}.

\smallskip
%%%%%
Thus our aim is now to show that $R(s)$ is constant. As a preliminary step in this direction, by a direct analysis of the function $S_F(s)$ in \eqref{5-14} we show that $R(s)$ is constant if $N=1$ and $N=2$; see Lemma 5.3 of \cite{Ka-Pe/2022b}. Note, in view of \eqref{5-16}, that $N\geq 3$ implies $\omega_{N-1}>0$. Hence we summarise the output of the second step by the following

\medskip
%%%%%
{\bf Proposition 4.1.} {\sl Let $F$ be as in Theorem $1$ and $\ga(s)$ be the virtual $\ga$-factor associated with $F$. Then $F$ satisfies functional equation \eqref{virtual-fe}, and Theorem $1$ follows if $R(s)$ is constant. Moreover, if $R(s)$ is not constant then, recalling definition \eqref{5-14}, $N\geq3$ and $\omega_{N-1}>0$.}

\medskip
%-5-%%%%%%%%%%%%%%%%%%%%%%%%%%%%%%%%%%%%%%%%%%%%%%%%%%%%%%%%%%%%%%%%%%%%%%%%%%%%%%%%%%%%%%%%%%%%
\section{Step 3: period functions}

\smallskip
%%%%%
In the last two steps we prove that the assumption $R(s)\neq$ constant leads to a contradiction, thus proving Theorem 1. This, however, requires several arguments. Our approach is based on a close study of the analytic properties of the Fourier series
\[
f(z) = \sum_{n=1}^\infty a(n) n^\lambda e(nz), \quad \text{where} \quad \lambda=
\begin{cases}
\mu & \text{if $\chi_F >0$ (Hecke case)} \\
i\kappa & \text{if $\chi_F \leq 0$ (Maass case)}
\end{cases}
\]
with $\mu$ and $\kappa$ as in \eqref{5-1}. Clearly, $f(z)$ is holomorphic on the upper half-plane $\HH = \{z\in\CC:\Im(z)>0\}$. Roughly, the function $f$ plays the role of the Hecke form associated to a function $F$ satisfying functional equation \eqref{virtual-fe} involving the virtual $\ga$-factor $\ga(s)$. Thus, according to Hecke's classical paradigm, we start expressing $f(iy)$ as a Mellin transform involving $F$ and then apply such a functional equation, trying to get a kind of modularity for $f$. More precisely, we shall eventually obtain a period function 
\begin{equation}
\label{perfun}
\psi(z) = f(z) - z^{-2\la-1} f(-1/z)
\end{equation}
associated with $f$ in the sense of Lewis-Zagier \cite{Le-Za/2001}, i.e. satisfying the three-term functional equation in \eqref{6-23} below, with $\mu=\lambda$. Moreover, we shall prove that the period function $\psi(z)$, clearly holomorphic on $\HH$, has continuation to $|\arg(z)|<\pi$.

\smallskip
%%%%%
Arguing in this way we first obtain that
\begin{equation}
\label{6-9}
f(iy) =  \frac{\omega_F}{2\pi i} \int_{(c_1)} (2\pi)^{-s} \Gamma(s) R(s-\lambda) \frac{\ga(1-s+\lambda)}{\ga(s-\lambda)} F(1-s+\lambda) y^{-s} \d s + L(iy),
\end{equation}
where $L(z)$ is a harmless residual term coming from the possible pole of $F$ at $s=1$ and is holomorphic for $|\arg(z)|<\pi$, and $c_1=\mu-\delta$ with a sufficiently small $\delta>0$ if $\la=\mu$ and $0<c_1<1$ if $\la=i\kappa$. From now on, the technical treatment of the Hecke and Maass cases is somewhat different. We discuss only the former case, i.e. $\chi_F>0$, the latter being simpler thanks to the shape of the corresponding virtual $\ga$-factor. Thus, assuming that $\la=\mu$, recalling the shape of $\ga(s)$ in \eqref{5-1} and observing that $S_\ga(s) = 2\sin(\pi(s+\mu))$ in this case, after some manipulations we transform \eqref{6-9} to
\begin{equation}
\label{6-11}
f(iy) = \omega_F (2\pi)^{1-\mu} y^{-\mu-1} \sum_{n=1}^\infty a(n) \frac{1}{2\pi i} \int_{(1+\delta)} \frac{1}{S_F(s)\Ga(1-s-\mu)}  \Big(\frac{2\pi n}{y}\Big)^{-s} \d s + L(iy).
\end{equation}

\smallskip
%%%%%
Next we observe that by \eqref{5-14} and \eqref{5-16}
\begin{equation}
\label{S-bound}
S_F(s) = -2\omega_F \cos(\pi s) + \sum_{j=1}^{N-1} a_j e^{i\pi\omega_js} \quad \text{and} \quad \frac{1}{S_F(s)} +\frac{1}{2\omega_F\cos(\pi s)} \ll e^{-\pi(2-\omega_{N-1})|t|}.
\end{equation}
Thus, isolating the main term $-2\omega_F \cos(\pi s)$, \eqref{6-11} can be rewritten as
\begin{equation}
\label{6-13}
\begin{split}
f(iy) &= -\frac12 (2\pi)^{1-\mu} y^{-\mu-1} \sum_{n=1}^\infty a(n) J\Big(\frac{2\pi n}{y}\Big) + H(iy) + L(iy) \\
&= \widetilde{K}(iy) + H(iy) +L(iy),
\end{split}
\end{equation}
say, where
\[
J(w) = \frac{1}{2\pi i} \int_{(1+\delta)} \frac{1}{\cos(\pi s)\Gamma(1-s-\mu)} w^{-s}\d s
\]
and
\begin{equation}
\label{6-5}
\begin{split}
H(z) = \omega_F(2\pi)^{1-\mu} (-iz)^{-\mu-1} \sum_{n=1}^\infty a(n) &\frac{1}{2\pi i} \int_{1+\delta} \Big(\frac{1}{S_F(s)}+\frac{1}{2\omega_F\cos(\pi s)}\Big) \\
&\times \frac{1}{\Ga(1-s-\mu)} \Big(\frac{2\pi in}{z}\Big)^{-s} \d s.
\end{split}
\end{equation}
Note that the function $H(z)$ is denoted by $Q_H(z)$ in \cite{Ka-Pe/2022b}. Note also that \eqref{S-bound} shows that $H(z)$ is holomorphic for 
\begin{equation}
\label{H-range}
-\pi(1-\omega_{N-1}) < \arg(z) < \pi,
\end{equation}
while $\widetilde{K}(z)$ is holomorphic for $z\in\HH$. 

\smallskip
%%%%%
At this point we anticipate that the argument we shall present in the next section leads to the conclusion of the proof of Theorem 1 provided we can show that $H(z)$ has holomorphic continuation to any angular region larger than \eqref{H-range}, i.e. of the form $-\pi\rho<\arg(z)<\pi$ with any fixed $\rho>1-\omega_{N-1}$. So, the remaining part of this section is devoted to this goal. Actually, we prove that $H(z)$ has continuation to $|\arg(z)|<\pi$; this is achieved by a detailed study of the function $\widetilde{K}(z)$, which brings into play the period function $\psi(z)$ and its properties.

\smallskip
We start the study of $\widetilde{K}(z)$ by rewriting the function $J(w)$ as a sum of residues shifting the line of integration to $-\infty$, namely
\begin{equation}
\label{6-14}
J(w) = -\frac{w^{-1/2}}{\pi} \sum_{\ell=0}^\infty \frac{(-w)^\ell}{\Gamma(\ell + 1/2 -\mu)} = -\frac{w^{-1/2}}{\pi} E_{1/2-\mu}(w),
\end{equation}
say. Here, in view of the relations between $w$ and $z$ implicit in \eqref{6-13}, we choose the branch of $w^{-1/2}$ with $-\pi/2 < \arg(w) < 3\pi/2$. Then it turns out that $E_{1/2-\mu}(w)$ is a classical special function, precisely a special case of two-parametric Mittag-Leffler function; see Chapter 4 of \cite{GKMR/2014}. Its relevant properties are that it is an entire function and has the more explicit representation
\begin{equation}
\label{6-15}
E_{1/2-\mu}(w) = \kappa_0 e^{-w} w^{1/2 +\mu} + \frac{e^{-w} w^{1/2 +\mu}}{\Ga(-1/2-\mu)}I_{1/2-\mu}(w),
\end{equation}
where $\kappa_0$ is a certain constant and  $I_{1/2-\mu}(w)$ is a certain absolutely convergent integral.

\smallskip
%%%%%
Now we plug \eqref{6-14} and \eqref{6-15} with $w=2\pi n/y$ into \eqref{6-13} and observe that the first term on the right hand side of \eqref{6-15} gives raise to $f(i/y)$, while the second one produces a term which we denote by $K(iy)$. Moreover, a direct estimate shows that $K(z)$ is holomorphic for $|\arg(z)|<\pi$. Thus, for $z\in\HH$ and writing $\kappa_1=\kappa_0 i^{2\mu+1}$, \eqref{6-13} transforms to
\begin{equation}
\label{6-18}
\begin{split}
f(z) &= \kappa_1 z^{-2\mu-1} f(-1/z) +K(z) + H(z) +L(z) \\
&= \kappa_1 z^{-2\mu-1} f(-1/z) +\widetilde{\psi}(z),
\end{split}
\end{equation}
say, where $\widetilde{\psi}(z)$ is holomorphic in the same range \eqref{H-range} as $H(z)$, since $K(z)$ and $L(z)$ are holomorphic for $|\arg(z)|<\pi$.

\smallskip
%%%%%
Next, by an argument based on the 1-periodicity of $f$, from \eqref{6-18} we deduce the relation
\begin{equation}
\label{6-20}
\begin{split}
 \widetilde{\psi}(z) - \frac{1}{\kappa_1} \widetilde{\psi}(z+1) &= \Big(1-\frac{1}{\kappa_1}\Big) f(z) - \kappa_1 z^{-2\mu-1} f(-1/z) + (z+1)^{-2\mu-1} f\Big(-\frac{1}{z+1}\Big) \\
 &= \Big(1-\frac{1}{\kappa_1}\Big) f(z) + (z+1)^{-2\mu-1} \widetilde{\psi}\Big(\frac{z}{z+1}\Big).
 \end{split}
\end{equation}
In turn, from \eqref{6-20} we deduce that $\kappa_1=1$. Indeed, all the $\widetilde{\psi}$-functions in \eqref{6-20} are holomorphic in an angular region of type $|\arg(z)|<\pi\rho$, thus if $\kappa_1\neq1$ then $f(z)$ would be entire thanks to its 1-periodicity. But this easily leads to a contradiction, so $\kappa_1=1$. Therefore, thanks to \eqref{perfun} and \eqref{6-18} we see that $\widetilde{\psi}(z)=\psi(z)$, \eqref{6-20} becomes the three-term functional equation
\begin{equation}
 \label{6-23}
\psi(z) = \psi(z+1) + (z+1)^{-2\mu-1} \psi\big(\frac{z}{z+1}\big)
\end{equation}
and $\psi(z)$ is holomorphic in the range \eqref{H-range}. Finally, we exploit \eqref{6-23} to show by an elementary geometric argument that $\psi(z)$ has holomorphic continuation to $|\arg(z)|<\pi$. As a consequence, thanks to \eqref{6-18}, $H(z)$ is also holomorphic for $|\arg(z)|<\pi$.

\smallskip
%%%%%
We summarise the output of the third step by the following proposition.

\medskip
%%%%
{\bf Proposition 5.1.} {\sl Let $F$ be as in Theorem $1$ and $\chi_F>0$. Then the function $H(z)$ in \eqref{6-5} is holomorphic for $|\arg(z)|<\pi$.}

\medskip
%%%%%
As we already pointed out, a similar result holds in the Maass case as well, i.e. if $\chi_F\leq0$.

\medskip
%-6-%%%%%%%%%%%%%%%%%%%%%%%%%%%%%%%%%%%%%%%%%%%%%%%%%%%%%%%%%%%%%%%%%%%%%%%%%%%%%%%%%%%%%%%%%%%%
\section{Step 4: conclusion of the proof}

\smallskip
%%%%%
Again, we discuss only the Hecke case, i.e. $\chi_F>0$, the other case being similar. Thanks to Proposition 4.1 we may assume that the function $R(s)$ is not constant, otherwise Theorem 1 follows, and hence that $N\geq3$ and $\omega_{N-1}>0$ in \eqref{5-14}. 

\smallskip
%%%%%
Our strategy is now as follows. The function $H(z)$ in \eqref{6-5} can be thought as a sum for $j=1,\dots,N-1$ of the terms corresponding to the frequencies $\omega_j$ in \eqref{S-bound}. Then we single out the term corresponding to $j=1$ and, recalling the symmetry of the $\omega_j$ in \eqref{5-16}, we show that it rebuilds (roughly) the function $f\left(\frac{-1}{e^{i\pi(1-\omega_{N-1})}z}\right)$. Moreover, a direct estimation shows that the sum of the remaining terms is holomorphic for $-\pi(\min(1-\omega_{N-2},1)) < \arg(z) < \pi$. But $H(z)$ is holomorphic for $|\arg(z)|<\pi$ and $\omega_{N-2}<\omega_{N-1}$, so by an obvious change of variable we see that $f(-1/z)$ is holomorphic in the region $-\pi\rho<\arg(z)<\pi$ for some $\rho>0$. By 1-periodicity this shows that $f(z)$ is entire, but this implies that actually $f(z)$ vanishes identically, a contradiction. Theorem 1 thus follows. 

\smallskip
%%%%%
Note, as already pointed out in the previous section, that what really matters on $H(z)$ is that it has holomorphic continuation to any angular region larger than \eqref{H-range}, which is enough to ensure that $f(-1/z)$ is holomorphic in a region of type $-\pi\rho<\arg(z)<\pi$ with $\rho>0$ and hence to get the same conclusions.

\smallskip
%%%%%
Some details of the above arguments are provided in the next section, which uses the notation in the original paper \cite{Ka-Pe/2022b}. However, the main difference in notation is that the function $H(z)$ of the present paper is denoted by $Q_H(z)$ in \cite{Ka-Pe/2022b}.

\medskip
%-7-%%%%%%%%%%%%%%%%%%%%%%%%%%%%%%%%%%%%%%%%%%%%%%%%%%%%%%%%%%%%%%%%%%%%%%%%%%%%%%%%%%%%%%%%%%%%
\section{Corrigendum of \cite{Ka-Pe/2022b}}

\smallskip
%%%%%
Due to an unfortunate slip in the definition of the function $G_H(z,s)$ in Section 7 of \cite{Ka-Pe/2022b}, the part where $\chi_F>0$ of that section needs to be corrected at several places. However, such a slip does not affect the basic lines of the argument, and the final result remains unchanged. 

\smallskip
%%%%%
Here we sketch the required changes in such a Section 7, from ``Suppose first that $\chi_F>0$'' to ``thus finishing the proof of the Hecke case.''. In this section we refer to notation, equation numbering and other results in \cite{Ka-Pe/2022b}.

\smallskip
%%%%%%%
The correct definition of $G_H(z,s)$ is
\[
G_H(z,s) = \ga(s) F(s) (-iz)^s.
\]
Then we split the integral in (7.1) exactly as in (7.2), but the role of $Q_H^-(z)$ and $Q_H^+(z)$ is now swapped, due to the change of sign of the exponent of $-iz$ in $G_H(z,s)$. So, we bound the integrand in $Q_H^-(z)$ by the standard estimate
\[
\ll (t+1)^c e^{(\arg(z)-\pi)|t|},
\]
thus getting that (7.3) holds with $Q_H^+(z)$ replaced by $Q_H^-(z)$.

\smallskip
Next we split $Q_H^+(z)$ similarly as in (7.4), taking into account that now integration is over $t\in[0,+\infty)$. So, recalling the symmetry of the $\omega_j$ in (5.15), equation (7.4) holds with $Q_H^-(z)$ replaced by $Q_H^+(z)$ and
\[
A_H(z) =  (2\pi)^{-\mu} (-iz)^{-\mu-1} \sum_{j=1}^{N-2} \frac{a_{N-j}}{2\pi i} \int_{1+\delta}^{1+\delta+i\infty} \frac{e^{-i\pi \omega_js}}{S_F(s)\cos(\pi s)} \sin(\pi(s+\mu)) G_H(z,s) \d s,
\]
\[
B_H(z) =  (2\pi)^{-\mu} (-iz)^{-\mu-1} \frac{a_1}{2\pi i} \int_{1+\delta}^{1+\delta+i\infty} \frac{e^{-i\pi \omega_{N-1}s}}{S_F(s)\cos(\pi s)} \sin(\pi(s+\mu)) G_H(z,s) \d s.
\]
Again, the integrand in $A_H(z)$ is
\[
\ll (|t|+1)^c e^{-|t|(\arg(z) + \pi -\pi\omega_j)},
\]
so $A_H(z)$ is holomorphic for $-\pi(\min(1-\omega_{N-2},1)) < \arg(z) < \pi$.

\smallskip
%%%%%
As in \cite{Ka-Pe/2022b}, the final step is the transformation of $B_H(z)$. Equation (7.6) now becomes
\[
\frac{\sin(\pi(s+\mu))}{S_F(s)\cos(\pi s)} =c_1 e^{i\pi s} \big( 1 + O(e^{-\pi(1-\omega_{N-1})|t|})\big) \qquad \text{as} \  t\to+\infty,
\]
hence we write
\[
B_H(z) = c_2 (-iz)^{-\mu-1} \frac{1}{2\pi i} \int_{1+\delta}^{1+\delta+i\infty} e^{i\pi s(1-\omega_{N-1})} G_H(z,s) \d s + C_H(z)
\]
with a certain $C_H(z)$ as in (7.8). Then, since as $t\to-\infty$ the integrand in the last integral is
\[
\ll (|t|+1)^c e^{|t|(\arg(z) -\pi\omega_{N-1})},
\]
the extension of such an integral to the whole line $\si=1+\delta$ introduces a further error term which, as in \cite{Ka-Pe/2022b}, is holomorphic for $-\pi<\arg(z)<\pi\omega_{N-1}$. Therefore, writing
\[
I_H(z) =  \frac{1}{2\pi i} \int_{(1+\delta)} e^{i\pi s(1-\omega_{N-1})} G_H(z,s) \d s
\]
and recalling that $Q_H(z)$ is holomorphic for $|\arg(z)|<\pi$ thanks to Proposition 6.1, gathering all the above results we obtain that $I_H(z)$ is holomorphic in the range (7.9). But in view of (6.8) we have that
\[
I_H(z) = (2\pi i)^\mu \big(e^{i\pi(\omega_{N-1}-1)}z\big)^{-\mu} f\left(\frac{-1}{e^{i\pi(1-\omega_{N-1})}z}\right).
\]
Hence, recalling (7.9) and that $\omega_{N-2}<\omega_{N-1}$, by the substitution $w=e^{i\pi(1-\omega_{N-1})}z$ we see that $f(-1/w)$ is holomorphic for $-\pi\rho<\arg(w)<\pi$ with some $0<\rho<1$. This, as in the paper, leads to a contradiction in view of Lemma 6.1, thanks to the 1-periodicity of $f$.

\medskip
%-8%%%%%%%%%%%%%%%%%%%%%%%%%%%%%%%%%%%%%%%%%%%%%%%%%%%%%%%%%%%%%%%%%%%%%%%%%%%%%%%%%%%%%%%%%%%%
\section{Proof of Theorem 2}

\smallskip
%%%%%
Here we show how Theorem 2 follows from Theorem 1. We denote by $a(n)$ and $b(n)$ the coefficients of $F$ and $G$, respectively. First we note that Theorem 2 reduces to Theorem 1 if $G=\pm \overline{F}$. Suppose now that $G\neq\pm \overline{F}$. We observe that every function $H_\pm(s) := c_\pm(F(s)\pm \overline{G}(s))$ with $c_\pm\in\CC\setminus\{0\}$ satisfies the functional equation
\[
\ga(s) H_\pm(s) = \omega_\pm \overline{\ga}(1-s) \overline{H_\pm}(1-s)
\]
with $\omega_\pm = \omega \overline{c_\pm}/c_\pm$, hence $H_\pm\in\S^\sharp$. Moreover, $d_{H_\pm}=d_\ga=2$, $q_{H_\pm}=q_\ga =1$ and $\theta_{H_\pm}=\theta_\ga$. Hence the functions $H_\pm$ can be normalized by means of a common shift by $-i\theta_\ga$ and choosing $c_\pm = (a(n_\pm) \pm \overline{b(n_\pm)})^{-1}$, where $n_\pm$ is the least value of $n$ such that $a(n)\pm \overline{b(n)}\neq 0$. We also recall that if one of the two normalized functions $H_\pm$  has a pole at $s=1$, then we already have that $\theta_\ga=0$; see Step 1 in the proof of Theorem 1. From now on we consider the normalized functions $H_\pm$, and note that $\chi_{H_\pm}=\chi_\ga$. Theorem 2 will follow, essentially, from an application of Theorem 1 to the functions $H_\pm$.

\smallskip
%%%%%
Assume that $\chi_\ga>0$ and let $=S_k(\Ga_0(1))$ denote the space of holomorphic cusp forms of level 1 and even integral weight $k= 1 + \sqrt{2\chi_\ga}$. Then, according to Theorem 1, there exist $f_+,f_-\in S_k(\Ga_0(1))$ such that
\begin{equation}
\label{T2-1}
H_\pm(s) = L\big(s+\frac{k-1}{2},f_\pm\big).
\end{equation}
Hence, denoting by $\overline{g}$ the cusp form having as coefficients the conjugates of those of $g$, solving the linear system \eqref{T2-1} in $F$ and $\overline{G}$ we see that there exist $f,\overline{g}\in S_k(\Ga_0(1))$ such that
\begin{equation}
\label{T2-2}
F(s) = L\big(s+\frac{k-1}{2},f\big), \quad \overline{G}(s) = L\big(s+\frac{k-1}{2},\overline{g}\big).
\end{equation}
But from Step 2 in the proof of Theorem 1 we know that $\ga(s) =(2\pi)^{-s} \Ga(s+\mu)$ with a certain real parameter $\mu$ depending on $\chi_\ga$ (actually, $\mu=\frac{k-1}{2}$). Hence, recalling the functional equation satisfied by the $L$-functions of holomorphic cusp forms, from \eqref{T2-2} and \eqref{intro-2} we obtain that
\[
\begin{split}
\ga(s) L\big(s+\frac{k-1}{2},f\big) &= \omega \ga(1-s) L\big(1-s +\frac{k-1}{2},g\big) \\
\ga(s) L\big(s+\frac{k-1}{2},f\big) &= i^k \ga(1-s) L\big(1-s +\frac{k-1}{2},f\big).
\end{split}
\]
Therefore $g= i^{k} \overline{\omega} f$, and the first claim of Theorem 2 follows by shifting back by $i\theta_\ga$.

\smallskip
Suppose now that $\chi_\ga=0$. By a similar argument we obtain that $F=\alpha\zeta^2$ and $G=\beta\zeta^2$ with some $\alpha,\beta\in\CC\setminus\{0\}$, and comparing with the functional equation of $\zeta^2$ we see that $\beta=\alpha\overline{\omega}$. Our second claim follows.

\smallskip
%%%%%
The third case, $\chi_\ga<0$, is also similar to the first one, with the difference that the root number of the functional equation in this case is $(-1)^\ep$, $\ep$ being the parity of the form. So in this case $g=(-1)^\ep \overline{\omega} f$, $\ep$ being as in \eqref{intro-0}, and the proof is complete. \qed

%\newpage
%-REFERENCES-%%%%%%%%%%%%%%%%%%%%%%%%%%%%%%%%%%%%%%%%%%%%%%%%%%%%%%%%%%%%%%%

\ifx\undefined\bysame{poly}.
\newcommand{\bysame}{\leavevmode\hbox to3em{\hrulefill}\ ,}
\fi

%\bigskip
%\bigskip
\medskip
\noindent
Jerzy Kaczorowski, Faculty of Mathematics and Computer Science, A.Mickiewicz University, 61-614 Pozna\'n. e-mail: \url{kjerzy@amu.edu.pl}

\medskip
\noindent
Alberto Perelli, Dipartimento di Matematica, Universit\`a di Genova, via Dodecaneso 35, 16146 Genova, Italy. e-mail: \url{perelli@dima.unige.it}

\end{document}